\documentclass[article,reqno]{amsart}
\usepackage[utf8]{inputenc}
\usepackage{a4wide, amsfonts}
\usepackage[english]{babel}
\usepackage{xcolor}
\usepackage[T1]{fontenc}
\usepackage{amsmath, amssymb,epic,graphicx,mathrsfs,enumerate}    
\usepackage{cite}
\usepackage[numbers,sort&compress]{natbib}
\usepackage[all]{xy}
\usepackage{color}
\usepackage{comment}
\usepackage{enumitem}
\usepackage{esint}
\usepackage{hyperref}
\usepackage{amsthm}
\usepackage{latexsym}
\usepackage{epsfig}
\usepackage{soul}
\usepackage{geometry}
\usepackage{doi}
\hypersetup{colorlinks=true, linkcolor=black, filecolor=black, urlcolor=blue,
citecolor=blue}

\newtheorem{theorem}{Theorem}[section]
\newtheorem{corollary}[theorem]{Corollary}

\newtheorem{prop}[theorem]{Proposition}

\theoremstyle{remark}
\newtheorem{obs}{Remark}

\newcommand{\R}{\mathbb{R}}

\newcommand{\e}{\varepsilon}
\newcommand{\al}{\alpha}

\newcommand{\s}{\sigma}

\title{Boundary value problems for Choquard equations}
\author[C.Bernardini - A.Cesaroni]{Chiara Bernardini and Annalisa Cesaroni}
\address{Dipartimento di Matematica "Tullio Levi Civita",
Università di Padova, Via Trieste 121, 35121 Padova, Italy}
\email{chiara.bernardini@math.unipd.it, annalisa.cesaroni@unipd.it}

\begin{document}

\begin{abstract}
\noindent We prove existence of a positive radial solution to the Choquard equation 
$$-\Delta u +V u=(I_\al\ast |u|^p)|u|^{p-2}u\qquad\text{in}\,\,\,\Omega$$ 
with Neumann or Dirichlet boundary conditions, when  $\Omega$ is an annulus, or an exterior domain of the form $\R^N\setminus \bar{B}_a(0)$. We provide also a nonexistence result, that is if $p\ge\frac{N+\al}{N-2}$ the corresponding Dirichlet problem does not have any nontrivial regular solution in strictly strictly star-shaped domains. \\
%Finally, we show that passing to the limit as $\al\to0^+$ we obtain an existence result for the corresponding \textit{local} problem with power-type nonlinearity.\\


\end{abstract}

\maketitle

\date{}

%-----------------------------------------------------------------------------

%\tableofcontents

\section{Introduction}

In the present work, we study the following nonlinear Choquard equation
\begin{equation}\label{choquard}
\begin{cases}
-\Delta u +V(x)u=(I_\al\ast |u|^p)|u|^{p-2}u \hspace{0.8cm}\text{in}\,\,\,\Omega\\
u>0\hspace{5.23cm}\text{in}\,\,\,\Omega
\end{cases}
\end{equation}
with Dirichlet or Neumann boundary conditions. We assume that $N\in\mathbb{N}$, $N\ge2$, the exponent in the nonlinearity is a real value $p>1$ and the potential $V:\Omega\to\R$ is a continuous radial function such that $\inf_{x\in\Omega}V>0$. Here, $I_\al:\R^N\to\R$ is the Riesz potential of order $\al\in(0,N)$, which is defined for every $x\in\R^N\setminus\{0\}$ by 
\begin{equation}\label{def riesz potential}
I_\al(x)=\frac{C_{N,\al}}{|x|^{N-\al}},\quad\text{where}\quad C_{N,\al}=\frac {\Gamma\left(\frac{N-\al}{2}\right)}{\Gamma(\frac{\al}{2})\pi^\frac{N}{2}2^\al},   
\end{equation}
recall that by $I_\al\ast f$ in a domain $\Omega\subset\R^N$, we mean the convolution 
$I_\al\ast(\chi_{\Omega}f)$ in $\R^N$. We consider both the case when the domain $\Omega$ is an annulus centered at the origin, namely   
$$\Omega=A_{a,b}:=\{x\in\R^N\,|\,a<|x|<b\,\}\qquad \text{ for some}\,\,\,0<a<b<+\infty$$ 
and the case of the exterior domain 
$$\Omega=\R^N\setminus\overline{B_a(0)}=\{x\in\R^N\,|\,|x|>a>0\,\}.$$ 

The Choquard equation has been extensively studied over the last decades, since it arises in the  modeling of several  mean-field physical phenomena. More in detail, the Choquard-Pekar equation 
\begin{equation}\label{cho-pek}
-\Delta u+u=(I_2\ast |u|^2)u \qquad\text{in}\,\,\,\R^3,    
\end{equation}
has been first introduced by S. Pekar \cite{Pekar} in 1954 to model the quantum mechanics of a polaron at rest, then P. Choquard used it to describe an electron trapped in its own hole. A further application was proposed by R. Penrose, who used it to model self-gravitating matter.\\
Existence of solutions to equation \eqref{cho-pek} in the normalized framework, namely imposing that $\|u\|_{L^2(\R^3)}=\mu$, has been first investigated using variational methods by E.H. Lieb \cite{Lieb} and in more general cases by P.-L. Lions \cite{Lions} (refer also to \cite{BerCes,Bernardini4} for further results about existence of solutions to more general Choquard-type systems). In particular, using symmetric decreasing rearrangement inequalities, E.H. Lieb proved that there exists a minimizing solution, which is radial and unique up to translations, while more recently, L. Ma and L. Zhao \cite{MaZhao} classified all positive solutions to \eqref{cho-pek}. 
On the other hand, the Choquard equation on $\R^N$ with a more general nonlocal nonlinearity depending on a parameter $p>1$, that is the following semilinear elliptic equation
\begin{equation}\label{cho p}
-\Delta u+u=(I_\al\ast |u|^p)|u|^{p-2}u \qquad\text{in}\,\,\,\R^N,    
\end{equation}
admits a nontrivial solution $u\in H^1(\R^N)\cap L^\frac{2Np}{N+\al}(\R^N)$ with $\nabla u\in H^1_{loc}(\R^N)\cap L^\frac{2Np}{N+\al}_{loc}(\R^N)$ if and only if $\frac{N+\al}{N}<p<\frac{N+\al}{N-2}$ (refer to \cite{MVS13,MS}  for a complete overview on the topic). The situation changes  when adding  an external potential $V$ in   the equation, that is considering
\begin{equation}\label{cho V}
-\Delta u+Vu=(I_\al\ast |u|^p)|u|^{p-2}u \qquad\text{in}\,\,\,\R^N,
\end{equation}
since  the presence of the (possibly) variable potential $V$  influences particles and hence it could affect existence of solutions. We refer the reader to  \cite{MVS15 2,VanSchXia} and \cite[§4]{MS}, for further discussion and references on the Choquard equation on $\R^N$ with a nonconstant potential $V$.
As we mentioned above, Choquard equations on the whole space $\R^N$ have been extensively studied in past decades, while there are few results about Choquard equations on other types of domain $\Omega\subset\R^N$. Notice that the notion of solution to the Choquard equation is \textit{nonlocal}, that is if $u\ge0$ (weakly) solves \eqref{cho V},  than $u$ is only a supersolution to the same equation in $\Omega\subset\R^N$. We refer, among others, to   \cite{GaoYang} for a Brezis-Nirenberg type critical problem of the Choquard equation on bounded domains, to  \cite{GhiPag} who proved that the number of positive solutions depends on the topology of the domain  for slightly subcritical Choquard problems, and to \cite{Goel1,Goel2} for existence results at the HLS critical level for  problems on non-contractible domains which contain a sufficiently large annulus. On the other hand, for what concerns exterior domains, we mention the work of V. Moroz and J. Van Schaftingen \cite{MVS13bis} regarding sharp Liouville-type nonexistence results for supersolutions in some suitable range of the parameter $p$, and optimal decay rates for solutions (see also \cite{ChenLiu,ClappSal} and the references therein).\\

In this paper, we prove that \eqref{choquard} admits a positive radial solution   both with Neumann and Dirichlet boundary conditions  in  annular domains $\Omega=A_{a,b}$ for every $p\geq 1$ and  in exterior domains $\Omega=\R^N\setminus \bar{B_a}$ for every $p\geq\frac{N+\alpha}{N}$, so generalizing to the case of nonlocal nonlinear equations  a classical result by Kadzan and Warner \cite{KadzanWarner} about existence of solutions to nonlinear equations with power-like nonlinearities in annular domains.  

Problem \eqref{choquard} has a variational structure: weak solutions are formally critical points (with $u>0$ on $\Omega$) of the action functional $\mathcal{A}$ defined for a function $u:\R^N\to\R$ by
\begin{equation}\label{funzionale azione A}
\mathcal{A}(u)=\frac{1}{2}\int\limits_{\Omega}\left(|\nabla u|^2+V(x) u^2\right)- \frac{1}{2p}\int\limits_{\Omega}\left(I_\al\ast |u|^p\right)|u|^p.
\end{equation}
Notice that due to the classical Sobolev Embedding Theorem $H^1(\R^N)\hookrightarrow L^q(\R^N)$ for $q\in[2,2^*]$, the above energy is well-defined and sufficiently differentiable on the Sobolev space $H^1(\R)$ if $\frac{N+\al}{N}\le p\le \frac{N+\al}{N-2}$. Setting the problem in an annular domain or in the exterior of a ball, and looking for radial solutions, we are able to enlarge the range of parameters $p$ for which the energy is well-defined. We will consider the Sobolev space of radial functions
$$H^1_{rad}(\Omega):=\{u\in H^1(\Omega)\,\,|\,\,u(x)=u(|x|)\,\}$$
with the usual norm 
$$\|u\|_{H^1}:=\left(\|u\|^2_{L^2}+\|\nabla u\|^2_{L^2}\right)^\frac{1}{2}=
\left(\int_\Omega |u|^2+|\nabla u|^2\right)^\frac{1}{2}.$$
We take advantage of the fact that the embedding
\begin{equation}
H^1_{rad}(\Omega)\hookrightarrow L^q(\Omega)
\end{equation}
is compact for every $q\ge1$ if $\Omega=A_{a,b}$ and for every $q>2$ if $\Omega=\{x\in\R^N\,|\,|x|>a>0\}$ (see Theorem \ref{teo appendix 1} and Theorem \ref{teo appendix 2} below). Then, in order to find solutions to \eqref{choquard} with homogeneous Neumann or Dirichlet boundary conditions we consider a constrained variational problem (see also \cite{MVS15 2} and \cite{VanSchXia}). For every $\al\in(0,N)$ fixed, we look for minimizers of the energy functional 
\begin{equation}\label{energia Q}
Q(u)=\frac{1}{2}\int\limits_{\Omega}|\nabla u|^2+V u^2
\end{equation}
over the constrained set 
\begin{equation}\label{vincolo M alfa}
M_\al:=\bigg\{u\in H^1_{rad}(\Omega)\,\,\bigg|\,\,\int_{\Omega}\left(I_\al\ast u^p\right)u^p=1 \bigg\} .  
\end{equation}
Notice that the minimizers $u_\al$ can be chosen non-negative (possibly taking $|u_\al|$ since it holds $\int_\Omega\big|\nabla|u_\al|\big|^2dx= \int_\Omega|\nabla u_\al|^2dx$). Hence, if   
\begin{equation}\label{J alfa}
Q(u_\al)= J_\al:=\inf\{Q(u)\,|\, u\in M_\al\},
\end{equation}
up to multiplication by a constant, the function $u_\al$ is  a positive groundstate of \eqref{choquard} with Neumann homogeneous boundary conditions.

Our main results are the following. 
\begin{theorem}\label{teo esist Neumann}
Let $N\ge2$, $\al\in(0,N)$, $\Omega=A_{a,b}$ and $V(|x|)$ be a continuous radial function on $\Omega$ such that $\inf_{x\in\Omega}V(x)>0$. Then, for every $p\in[1,+\infty)$ there exists $u\in H^1_{rad}(\Omega)$ which solves 
\begin{equation}\label{choquad neumann}
\begin{cases}
-\Delta u +V(x)u=(I_\al\ast |u|^{p})|u|^{p-2}u \qquad\text{in}\,\,\,\Omega\\
u>0\hspace{5.12cm}\text{in}\,\,\,\Omega\\
\frac{\partial u}{\partial\nu}=0\hspace{4.87cm}\text{on}\,\,\,\partial\Omega
\end{cases}.
\end{equation}
On the other hand, if $\Omega=\{x\in\R^N\,|\,|x|>a>0\}$ the boundary value problem \eqref{choquad neumann} admits a solution for every $p\in\left(\frac{N+\al}{N},+\infty\right)$.
\end{theorem}

The analogous result holds in the case of Dirichlet boundary conditions. Notice that if $N\ge3$ and $\Omega$ is an annular domain, we can weaken the assumptions on the potential $V$, covering in this way also the case of the Choquard equation with the unperturbed Laplacian, that is when $V\equiv0$.

\begin{theorem}\label{teo esist Dirichlet}
Let $N\ge2$, $\al\in(0,N)$, $\Omega=A_{a,b}$ and $V(|x|)$ be a continuous radial function on $\Omega$ such that $\inf_{x\in\Omega}V(x)>0$ if $N=2$ and   $V(|x|)\ge0$ if $N\ge3$. Then, for every $p\in[1,+\infty)$ there exists $u\in H^1_{0,rad}(\Omega)$ which solves 
\begin{equation}\label{choquad dirichlet}
\begin{cases}
-\Delta u +V u=(I_\al\ast |u|^{p})|u|^{p-2}u \qquad\text{in}\,\,\,\Omega\\
u>0\hspace{4.65cm}\text{in}\,\,\,\Omega\\
u=0\hspace{4.62cm}\text{on}\,\,\,\partial\Omega
\end{cases}.
\end{equation}
If $\Omega=\{x\in\R^N\,|\,|x|>a>0\}$ the boundary value problem \eqref{choquad dirichlet} admits a solution when $\inf_{x\in\Omega}V(x)>0$ for every $p\in\left(\frac{N+\al}{N},+\infty\right)$.
\end{theorem}

On the other hand, when $\Omega$ is a smooth domain in $\R^N$ and the solution $u$ is sufficiently regular, we are able to obtain a Pohozaev-type identity, which in turns implies the triviality of $u$, when the domain $\Omega$ is strictly star-shaped.

\begin{theorem}\label{teo non esist}
Let $N\ge3$, $\Omega\subset\R^N$ be a smooth domain strictly star-shaped with respect to the origin, $p\ge\frac{N+\al}{N-2}$ and $V\in C^1(\Omega)$. Let $u\in H^1_0(\Omega)\,\cap H^2(\Omega)\,\cap W^{1,\frac{2pN}{N+\al}}(\Omega)$  be a solution to 
\begin{equation}\label{pb non exist}
\begin{cases}
-\Delta u +Vu=(I_\al\ast |u|^p)|u|^{p-2}u \qquad\text{in}\,\,\,\Omega\\
u=0\qquad\qquad\qquad\quad\,\qquad\qquad\,\,\,\,\,\quad\text{on}\,\,\partial\Omega.
\end{cases}
\end{equation}
Assume that 
\begin{equation}\label{cond 1 V}
u^2 V,\, u^2\nabla V\cdot x\in L^1(\Omega)   
\end{equation}
and 
\begin{equation}\label{cond 2 V}
\int_{\Omega}u^2\nabla V\cdot x\ge0.
\end{equation}
If $p=\frac{N+\al}{N-2}$, we assume that $V>0$ in $\Omega$, whereas if $p>\frac{N+\al}{N-2}$ it is sufficient that $V\geq 0$. 
Then $u\equiv 0$.
\end{theorem}

Finally, in the case $\Omega=A_{a,b}$, we prove a $\Gamma$-convergence type result as $\al\to 0^+$, relating the minimization problem \eqref{J alfa} with
$$J_0=\inf\{Q(u)\,|\,u\in M_0\}$$
where 
$$M_0:=\bigg\{u\in H^1_{rad}(\Omega)\,\,\bigg|\,\,\|u\|_{L^{2p}(\Omega)}=1,\quad u\geq 0\bigg\}.$$
This argument allows us to recover the existence result of Kadzan and Warner \cite{KadzanWarner} for solutions to the corresponding local problem
\begin{equation}\label{pb u2p}
\begin{cases}
-\Delta u+V u=u^{2p-1}\hspace{0.7cm}\text{in}\,\,\,\Omega\\
u>0\qquad\qquad\quad\quad\qquad\,\text{in}\,\,\,\Omega\\
\frac{\partial u}{\partial\nu}=0\hspace{2.3cm}\quad\text{on}\,\,\,\partial\Omega
\end{cases}
\end{equation}
and to extend it also to the case $V$ is nonconstant. Notice that, the analogous result holds also with Dirichlet boundary conditions, in this context we may refer among others to \cite{Grossi,NiNuss} for some related results.

\begin{theorem}\label{teo problema limite}
Let us assume that $\Omega=A_{a,b}$, $N\ge2$, $V(|x|)$ be a continuous radial function on $\Omega$ such that $\inf_{x\in\Omega}V(x)>0$ and $p$ is a fixed real value in $[1,+\infty)$. Then, 
$$\lim\limits_{\al\to0^+}J_\al=J_0$$
and if $\{u_\al\}$ is a sequence of minimizers for $J_\al$, there exists $u_0\in M_0$ such that as $\al\to0^+$ 
$$ u_\al\to u_0\,\,\,\text{in} \,\,\,H^1_{rad}(\Omega)$$ 
and $J_0=Q(u_0)$. Moreover, $u_0$ is a solution to \eqref{pb u2p} up to multiplication by a constant.
\end{theorem}

\begin{theorem}
Under the assumptions of Theorem \ref{teo esist Dirichlet}, the same result as Theorem \ref{teo problema limite} holds in the case of Dirichlet boundary condition on $\Omega=A_{a,b}$.
\end{theorem}

%-----------------------------------------------------------------------------------

\section{Sobolev embeddings for radial functions}

We state here some results on Sobolev embeddings for $H^1_{rad}(\Omega)$ and $H^1_{0,rad}(\Omega)$, both in the cases when the domain $\Omega\subset\R^N$ is an annulus  and when $\Omega=\{x\in\R^N\,|\,|x|>a>0\}$. Since we have not found a detailed proof in the literature, we recall it here by completeness.  

\begin{theorem}\label{teo appendix 1}
Let $N\ge2$ and $\Omega=A_{a,b}$. For every $p\in[1,+\infty)$, the following immersion 
$$H^{1}_{rad}(\Omega)\hookrightarrow L^p(\Omega)$$
is compact. Moreover, for every $p\ge1$ there exists a positive constant $C=C(N,p,a,b)$ such that 
\begin{equation}\label{stima sob anello}
\|u\|_{L^p(\Omega)}\le C\|u\|_{H^1_{rad}(\Omega)}.    
\end{equation}
\end{theorem}

\begin{obs}
If $a=0$ namely $\Omega=A_{0,b}=B_b(0)$, from the classical Rellich-Kondrachov Theorem (refer to \cite[Theorem 6.2]{Adams}) we have that the previous result holds for every $p\in[1,2^*)$ if $N\ge3$, and for every $p\in[1,+\infty)$ if $N=2$.
\end{obs}

\proof
If $N=2$, since $\Omega$ is bounded and smooth, we have that for all $p\in[1,+\infty)$ the embedding 
$$H^1(\Omega)\to L^p(\Omega)$$
is compact (see e.g. \cite[Theorem 6.2]{Adams}).

If $N\ge3$,   every function $u\in H^1_{rad}(A_{a,b})$ can be extended to a function $\bar{u}\in H^1_{rad}(\R^N)$ such that $\bar{u}|_{A_{a,b}}=u$ and $\bar{u}(x)=0$ for $|x|$ sufficiently large with 
%(by construction we have that $\bar{u}\equiv0$ for $|x|\ge2b-a$) 
$$\|\bar{u}\|_{L^2(\R^N)}\le C\|u\|_{L^2(A_{a,b})}\qquad\text{and}\qquad \|\bar{u}\|_{H^1(\R^N)} \le C\|u\|_{H^1(A_{a,b})}$$
where $C$ depends only on $|b-a|$,  see \cite[Theorem 8.6]{Brezis}. We recall also that every function $\bar{u}\in H^1(\R)$ is represented by a continuous function on $\Bar{\R}$, which we denote again by $\bar{u}$, and such that
$$\bar{u}(x)-\bar{u}(y)=\int_y^x\bar{u}'(r)\,dr,\qquad\forall x,y\in[-\infty,+\infty].$$
Hence following a classical computation by Strauss \cite{Strauss}, for any $u\in H^1_{rad}(A_{a,b})$ and $x\in A_{a,b}$ we have that
\begin{align*}
|u&(x)|^2=|\bar{u}(|x|)|^2=\bigg|\int_{|x|}^{+\infty}\frac{d}{dr}|\bar{u}(r)|^2 dr\bigg|\le\int_{|x|}^{+\infty}\bigg|\frac{d}{dr}|\bar{u}(r)|^2\bigg|dr\\
&\le2\left(\int_{|x|}^{+\infty}|\bar{u}(r)|^\frac{2N}{N-2}r^{N-1}dr\right)^\frac{N-2}{2N}\left(\int_{|x|}^{+\infty}\left|\frac{d}{dr}\bar{u}(r)\right|^2 r^{N-1}dr\right)^\frac{1}{2}\left(\int_{|x|}^{+\infty}r^{(2-N)N-1}dr\right)^\frac{1}{N}\\
&\le C|x|^{2-N}\|\bar{u}\|_{L^{2^*}(\R^N)}\|\nabla\bar{u}\|_{L^2(\R^N)}\le C|x|^{2-N}\|\nabla\bar{u}\|^2_{L^2(\R^N)}\\
&\le C |x|^{2-N}\|u\|^2_{H^1(A_{a,b})}
\end{align*}
where $C=C(N,a,b)$ and we used the H\"older's inequality and the classical Sobolev Embedding Theorem. Hence for every $p\in[1,+\infty)$ we get
$$\int_{A_{a,b}}|u(x)|^pdx=C\int_a^b r^{N-1}|u(r)|^p dr\le C \|u\|^p_{H^1(A_{a,b})} \int_a^b r^{N-1+p-\frac{Np}{2}}dr$$
which proves estimate \eqref{stima sob anello}. In order to prove compactness, let 
$\{u_n\,|\,n\in\mathbb{N}\}\subset H^1_{rad}(\Omega)$ be a bounded sequence, then there exists a constant $M>0$ such that $\|u_n\|_{H^1(\Omega)}\le M$ $\forall n\in\mathbb{N}$. Using \eqref{stima sob anello}, we get that the family $\{u_n\,|\,n\in\mathbb{N}\}$ is bounded in $L^p(\Omega)$ for any $p\in[1,+\infty)$. Let us show that $\{u_n\,|\,n\in\mathbb{N}\}$ is uniformly equi-continuous in $L^p(\Omega)$. For every $k\in\mathbb{N}_*$ we define
$$\Omega_k:=\left\{x\in\Omega\,\bigg|\,\mathrm{dist}(x,\partial\Omega)>\frac{1}{k}\right\}.$$
Using the H\"older's inequality and the Sobolev embedding, we can prove that $\forall k\in\mathbb{N}_*$ and $\forall h\in\R^N$ it holds
\begin{equation}\label{dim equicont 1}
\|\tau_hu-u\|_{L^1(\Omega\setminus\Omega_k)}\le C\|u\|_{H^1(\Omega)}\mathcal{L} (\Omega\setminus\Omega_k)^\frac{1}{2},   
\end{equation}
where $C=C(N,a,b)$ and $\tau_h u(x):=u(x+h)$ with $u$ extended with $0$ outside $\Omega$. Moreover $\forall k\in\mathbb{N}_*$ and $\forall h\in\R^N$ such that $|h|\le\frac{1}{k}$ we have
\begin{equation}\label{dim equicont 2}
\|\tau_hu-u\|_{L^1(\Omega_k)}\le|h| \, \|\nabla u\|_{L^2(\Omega)}\mathcal{L}(\Omega)^\frac{1}{2}.  
\end{equation}
Hence from \eqref{dim equicont 1} and \eqref{dim equicont 2} we get that $\forall k\in\mathbb{N}_*$, $\forall h\in\R^N$ such that $|h|\le\frac{1}{k}$ and $\forall n\in\mathbb{N}$ we have
\begin{align*}
\|\tau_hu_n-u_n\|_{L^1(\Omega)}&=\|\tau_hu_n-u_n\|_{L^1(\Omega\setminus\Omega_k)}+
\|\tau_hu_n-u_n\|_{L^1(\Omega_k)}\\
&\le M\left(C\mathcal{L}(\Omega\setminus\Omega_k) ^{1/2}+|h|\,\mathcal{L}(\Omega)^\frac{1}{2}\right).  
\end{align*}
For every fixed $\e>0$ we can choose $k_0$ such that if $|h|\le\frac{1}{k_0}$ then $$|h| \,\mathcal{L}(\Omega)^\frac{1}{2}\le\frac{\e}{M}\qquad\text{and}\qquad \mathcal{L}(\Omega\setminus\Omega_k)^{1/2}\le\frac{\e}{CM},$$
hence we get that $\forall\e>0$, $\exists\delta>0$ such that $\forall h\in\R^N$ with $|h|\le\delta$ and for any $n\in\mathbb{N}$ it holds
$$\|\tau_h u_n-u_n\|_{L^1(\Omega)}\le\e.$$
Finally for every $p\in(1,+\infty)$ fixed, by interpolation there exist $q\in(p,+\infty)$ and $\theta\in(0,1)$ such that 
\begin{align*}
\|\tau_h u_n-u_n\|_{L^p(\Omega)}&\le\|\tau_h u_n-u_n\|_{L^1(\Omega)}^\theta\|\tau_h u_n-u_n\|^{1-\theta}_{L^q(\Omega)}\\
&\le\|\tau_h u_n-u_n\|_{L^1(\Omega)}^\theta\left[\|\tau_h u_n\|_{L^q(\Omega)}+\|u_n \|_{L^q(\Omega)}\right]^{1-\theta}\\
&\le\|\tau_h u_n-u_n\|_{L^1(\Omega)}^\theta\left[2\|u_n \|_{L^q(\Omega)}\right]^{1-\theta}\\
&\le \|\tau_h u_n-u_n\|_ {L^1(\Omega)}^\theta (2CM)^{1-\theta} 
\end{align*}
where in the last inequality we used \eqref{stima sob anello}. This proves equicontinuity in $L^p(\Omega)$, namely $\forall\e>0$, $\exists\delta>0$ such that $\forall h\in\R^N$ with $|h|\le\delta$ and for any $n\in\mathbb{N}$ it holds
$$\|\tau_h u_n-u_n\|_{L^p(\Omega)}\le\e.$$ 
Then by the compactness criterion in $L^p$, we can conclude that the sequence $\{u_n\,|\,n\in\mathbb{N}\}$ converges (up to subsequences) in $L^p(\Omega)$.
\endproof

\begin{theorem}\label{teo appendix 2}
Let $N\ge2$ and $\Omega=\{x\in\R^N\,|\,\,|x|>a>0\}$. For every $p\in[2,+\infty)$ we have the following continuous immersion 
$$H^1_{rad}(\Omega)\hookrightarrow L^p(\Omega);$$
which is also compact for $p\in(2,+\infty)$.
\end{theorem}

\proof
\textit{Case $N\ge3$.} Let $u\in H^1_{rad}(\Omega)$, since it is radial we can identify it with a function in $H^1\big((a,+\infty)\big)$ which we still denote by $u$.
Since $u\in H^1((a,+\infty))$ then $\lim\limits_{|x|\to+\infty}u(x)=0$, by the same argument as in the proof of Theorem \ref{teo appendix 1}, we get that
\begin{equation}\label{stima u rad}
|u(x)|\le C|x|^{1-\frac{N}{2}}\|u\|_{H^1(\Omega)}
\end{equation}
where $C=C(N,a)$. If $p>2^*$, $N+p-\frac{Np}{2}<0$ and so
\begin{equation}\label{sob emb fuori palla p>2*}
\int_{\Omega}|u(x)|^p dx\le C_1\|u\|^p_{H^1(\Omega)}\int_a^{+\infty} r^{N-1+p- \frac{Np}{2}}dr=C_2\|u\|^p_{H^1(\Omega)}.   
\end{equation}
From \eqref{sob emb fuori palla p>2*} and the classical Sobolev embedding (see \cite[Theorem 5.4]{Adams}) which holds for $p\in[2,2^*]$, we can conclude that if $u\in H^1_{rad}(\Omega)$
$$\|u\|_{L^p(\Omega)}\le C_{N,p,a}\|u\|_{H^1(\Omega)},\qquad \forall p\in [2,+\infty)$$
and hence 
$$H^1_{rad}(\Omega)\hookrightarrow L^p(\Omega) ,\qquad \forall p\in [2,+\infty).$$
Now in order to prove compactness, we proceed as in \cite[Theorem 11.2]{AmbrosettiMalchiodi} (see also \cite{Strauss}). Let $\{u_n\,|\,n\in\mathbb{N}\}\subset H^1_{rad}(\Omega)$ be a bounded sequence, without loss of generality we can assume that $u_n\rightharpoonup0$ in $H^1_{rad}(\Omega)$. From \eqref{stima u rad} it follows that
$$|u_n(x)|\le C_1|x|^{1-\frac{N}{2}},\qquad\text{for every }n\in\mathbb{N}$$
so if $p>2$, given $\e>0$ there exist $C_2,R>0$ (we can always assume $R>a$) such that
$$|u_n(x)|^{p-2}\le C_1|x|^{\left(1-\frac{N}{2}\right)(p-2)}\le C_2\e,\qquad \text{for }|x|\ge R.$$
Using this we get that
\begin{equation}\label{magg R}
\int_{|x|\ge R}|u_n(x)|^p dx\le C_2 \e\int_{|x|\ge R}|u_n(x)|^2 dx\le C_2\e \|u_n\|^2_{H^1(\Omega)}\le C_3\e.    
\end{equation}
We consider now the annulus $A:=\{x\in\R^N\,|\, a<|x|<R\}$, since the sequence $\{u_n\,|\,n\in\mathbb{N}\}$ is bounded in $H^1_{rad}(A)$, recalling the compact embedding $H^1_{rad}(A)\subset\subset L^p(A)$ for any $p\in[1,+\infty)$ (see Theorem \ref{teo appendix 1}) we get that $u_n\to 0$ strongly in $L^p(A)$. It follows that there exists $n_0\in \mathbb{N}$ such that for all $n\ge n_0$
\begin{equation}\label{tra a e R}
\int_{a<|x|<R}|u_n(x)|^p dx\le \e.  
\end{equation}
Combining \eqref{magg R} and \eqref{tra a e R}, we get that for every $n\ge n_0$
$$\int_{\Omega}|u_n(x)|^p dx\le C_4\e$$
which proves that $u_n\to 0$ in $L^p(\Omega)$ $\forall p\in(2,+\infty)$.

\textit{Case $N=2$.} The continuous immersion $H^1(\Omega)\hookrightarrow L^p(\Omega)$ holds for every $p\in[2,+\infty)$ (see \cite[Theorem 5.4 Case B]{Adams}). For what concern compactness, we proceed as in the previous case using that for any $u\in H^1_{rad}(\Omega)$ it holds
$$|u(x)|\le C|x|^{-\frac{1}{4}}\|u\|_{H^1(\Omega)}.$$
Actually, if $N=2$ the embedding is compact for every $p\in[2,+\infty)$.
\endproof

\begin{obs}
Notice that if $N\ge3$ the embedding of $H^1_{rad}(\R^N)$ in $L^p(\R^N)$ is compact for all $p\in(2,2^*)$ (see \cite[Theorem 11.2]{AmbrosettiMalchiodi} and also \cite{Strauss}).
\end{obs}

Regarding the Sobolev space $H^1_{0,rad}(\Omega)$ we get the following results. 

\begin{corollary}\label{cor emb H_0 anello}
Let $N\ge2$ and $\Omega=A_{a,b}$. For every $p\in[1,+\infty)$, the following compact immersion holds
$$H^1_{0,rad}(\Omega)\hookrightarrow L^p(\Omega).$$
Moreover, if $N\ge3$ there exists a positive constant $C=C(N,p,a,b)$ such that $\forall p\ge1$
\begin{equation}\label{stima solo grad anello}
\|u\|_{L^p(\Omega)}\le C\|\nabla u\|_{L^2(\Omega)}.   
\end{equation}
\end{corollary}

\proof
It follows similarly to the proof of Theorem \ref{teo appendix 1}. If $N=2$, since $A_{a,b}$ is bounded, we have that for all $p\in[1,+\infty)$ the embedding 
$$H_0^1(\Omega)\hookrightarrow L^p(\Omega)$$
is compact (see e.g. \cite[Theorem 6.2]{Adams} which holds also in $H_0^1(\Omega)$).
If $N\ge3$ and $U\subset\R^N$ is an arbitrary open set, then for every $u\in H_0^1(U)$ it holds
\begin{equation}\label{SGN H_0^1}
\|u\|_{L^{2^*}(U)}\le C_N\|\nabla u\|_{L^2(U)} 
\end{equation} 
(classical Gagliardo-Nirenberg-Sobolev inequality for $H^1_0(U)$). Since we are working with functions in $H^1_{0,rad}(\Omega)$, extending $u$ by $0$ outside $\Omega$ and proceeding as before we obtain the following estimate
$$|u(x)|^2\le C|x|^{2-N}\|u\|_{L^{2^*}(\Omega)}\|\nabla u\|_{L^2(\Omega)}\le C|x|^{2-N}\|\nabla u\|^2_{L^2(\Omega)}$$
where in the last inequality we used \eqref{SGN H_0^1}. Estimate \eqref{stima solo grad anello} follows immediately; for what concerns compactness we can replicate the arguments before.
\endproof

\begin{obs}
Notice that differently from Theorem \ref{teo appendix 1}, if $N\ge3$ only the gradient of $u$ appears in the right-hand-side of inequality \eqref{stima solo grad anello}. This fact will be useful in the next section, in order to relax the assumptions on the potential $V$ for problems defined on annular domains.
\end{obs}

\begin{corollary}\label{cor emb H_0 fuoripalla}
Let $\Omega=\{x\in\R^N\,|\,|x|>a>0\}$. For every $p\in[2,+\infty)$ we have the following continuous immersion
$$H^1_{0,rad}(\Omega)\hookrightarrow L^p(\Omega).$$
Hence for every $p\ge2$ there exists a positive constant $C=C(N,p,a)$ such that 
\begin{equation}\label{stima sob0 fuori palla}
\|u\|_{L^p(\Omega)}\le C\|u\|_{H^1_{rad}(\Omega)}
\end{equation}
moreover, if $N\ge3$ for any $p\ge2^*$ we have
\begin{equation}\label{stima solo grad fuoripalla}
\|u\|_{L^p(\Omega)}\le C\|\nabla u\|_{L^2(\Omega)}.   
\end{equation}
The previous immersion is compact for every $p\in(2,+\infty)$.
\end{corollary}

\proof
\textit{Case $N\ge3$.} Let $u\in H^1_{0,rad}(\Omega)$, by the same argument as before and using estimate \eqref{SGN H_0^1}, we get that
\[
|u(x)|\le C|x|^{1-\frac{N}{2}}\|\nabla u\|_{L^2(\Omega)}
\]
where $C=C(N,a)$. If $p>2^*$, $N+p-\frac{Np}{2}<0$ and so
\[
\int_{\Omega}|u(x)|^p dx\le C_1\|\nabla u\|^p_{L^2(\Omega)}\int_a^{+\infty} r^{N-1+p- \frac{Np}{2}}dr=C_2\|\nabla u\|^p_{L^2(\Omega)}, 
\]
which proves estimate \eqref{stima solo grad fuoripalla}. By the previous inequality and the classical Sobolev embedding (see \cite[Theorem 5.4]{Adams}) which holds for $p\in[2,2^*]$, we can conclude that if $u\in H^1_{0,rad}(\Omega)$
$$\|u\|_{L^p(\Omega)}\le C_{N,p,a}\|u\|_{H^1(\Omega)},\qquad \forall p\in [2,+\infty)$$
and hence the continuous embedding follows. \\
\indent \textit{Case $N=2$.} The continuous immersion $H_0^1(\Omega)\hookrightarrow L^p(\Omega)$ holds for every $p\in[2,+\infty)$ (see \cite[Theorem 5.4 Case B]{Adams}).\\
\indent In both cases, the proof of compactness follows the same arguments as in the proof of Theorem \ref{teo appendix 2}.
\endproof

%------------------------------------------------------------------------------------

\section{Existence of a constrained minimizer}

We construct a radial solution to the Neumann boundary value problem \eqref{choquad neumann} as minimizer of $Q(u)$ on the constrained set $M_\al$ (as defined in \eqref{energia Q} and \eqref{vincolo M alfa}).

\begin{prop}\label{esist minimo anello neumann}
Let $N\ge2$, $\al\in(0,N)$ fixed, $\Omega=A_{a,b}$ and $V(|x|)$ be a continuous radial function on $\Omega$ such that $\inf_{x\in\Omega}V(x)>0$. Then, for every $p\in[1,+\infty)$, there exists $u_\al\in M_\al$ non-negative function such that $$Q(u_\al)=J_\al:=\inf\limits_{u\in M_\al}Q(u).$$
\end{prop}

\proof
Let $J_\al:=\inf\limits_{u\in M_\al} Q(u)\ge0$ and $\{u_n\,|\,n\in\mathbb{N}\} \in M_\al$ be a minimizing sequence for $J_\al$. We can assume that there exists $n_0\in\mathbb{N}$ such that for any $n\ge n_0$
$$Q(u_n)<J_\al+1$$
from which, using that $\inf\limits_{x\in\Omega}V(x)>0$, we deduce that
$$\int_{\Omega}|\nabla u_n|^2 dx<J_\al+1, \qquad \int_\Omega u_n^2\, dx\le \frac{Q(u_n)}{\inf V}<\frac{J_\al+1}{\inf V}.$$
This proves that the sequence $\{u_n\,|\,n\in\mathbb{N}\}$ is bounded in $H^1_{rad}(\Omega)$, hence there exists $u_\al\in H^1_{rad}(\Omega)$ such that $u_n\rightharpoonup u_\al$ weakly in $H^1_{rad}(\Omega) $ (up to subsequences) as $n\to\infty$ and almost everywhere in $\Omega$. Now, in order to conclude, we have to verify that $u_\al\in M_\al$. Using the Hardy-Littlewood-Sobolev inequality   we get that
\begin{align*}
\Bigg|\int_{\Omega}(I_\al&\ast|u_n|^p)|u_n|^p-\int_{\Omega}(I_\al\ast|u_\al|^p) |u_\al|^p\Bigg|\\
&=C_{N,\al}\Bigg|\,\int\limits_\Omega\int\limits_\Omega\frac{|u_n(x)|^p|u_n(y)|^p}{|x-y|^{N-\al}}dx\,dy-\int\limits_\Omega\int\limits_\Omega \frac{|u_\al(x)|^p |u_\al(y)|^p}{|x-y|^{N-\al}}dx\,dy\Bigg|\\
&=C_{N,\al}\Bigg|\,\int\limits_\Omega\int\limits_\Omega\frac{(|u_n(x)|^p-|u_\al(x)|^p)(|u_n(y)|^p+|u_\al(y)|^p)}{|x-y|^{N-\al}}dx\,dy\Bigg|\\
&\le C\Big\|\,|u_n|^p-|u_\al|^p\Big\|_{L^\frac{2N}{N+\al}(\Omega)}\Big\|\,|u_n|^p +|u_\al|^p\Big\|_ {L^\frac{2N}{N+\al}(\Omega)}\\
&\le C\Big\| |u_n|^p-|u_\al|^p\Big\|_{L^\frac{2N}{N+\al}(\Omega)}\left(\|u_n\|^p_ {L^\frac{2Np}{N+\al}(\Omega)}+ \|u_\al\|^p_{L^\frac{2Np}{N+\al}(\Omega)}\right).
\end{align*}
Using estimate \eqref{stima sob anello} and the fact that 
$\|u_n\|_{H^1(\Omega)}\le C$, we obtain that for every $p\in[1,+\infty)$
$$\|u_n\|^p_ {L^\frac{2Np}{N+\al}(\Omega)}+\|u_\al\|^p_{L^\frac{2Np}{N+\al}(\Omega)}\le C$$
uniformly in $n$. If $p=1$ we are done, since using Theorem \ref{teo appendix 1} we have that $\|u_n-u_\al\|_\frac{2N}{N+\al}\to 0$. On the other hand if $p>1$, in order to deal with the term $\||u_n|^p-|u_\al|^p\|_{L^ \frac{2N}{N+\al}(\Omega)}$, we take advantage of the following estimate. Using convexity of the function $x^p$, for $x>0$ and $p>1$, we get that
$$\Big| |u_n|^p-|u_\al|^p\Big|\le p\,\Big||u_n|-|u_\al|\Big|\,\Big(|u_n|^{p-1}+|u_\al|^{p-1}\Big)$$
and so using H\"older's inequality
\begin{align}\notag
\Big\| |u_n|^p-&|u_\al|^p\Big\|_{L^\frac{2N}{N+\al}(\Omega)}\\\notag
&\le C\Big\||u_n|-|u_\al|\Big\|_{L^\frac{2Nr}{N+\al}(\Omega)}\Big\| |u_n|^{p-1}+|u_\al|^{p-1}
\Big\|_{L^\frac{2Nr'}{N+\al}(\Omega)}\\ \label{formula stima diff ap-bp}
&\le C\Big\||u_n|-|u_\al|\Big\|_{L^\frac{2Nr}{N+\al}(\Omega)}\Big(\|u_n\|^{p-1} _{L^\frac{2Nr'(p-1)}{N+\al}(\Omega)}+\|u_\al\|^{p-1}_{L^\frac{2Nr'(p-1)}{N+\al}(\Omega)}\Big)
\end{align}
where $r$ and $r'$ are conjugate exponents. Choosing $r$ such that 
$$\frac{2Nr'(p-1)}{N+\al}\ge1,$$
we can use estimate \eqref{stima sob anello} again, and since $\{u_n\}$ is bounded in $H^1_{rad}(\Omega)$, we get that 
$$\|u_n\|^{p-1} _{L^\frac{2Nr'(p-1)}{N+\al}(\Omega)}+\|u_\al\|^{p-1}_{L^\frac{2Nr' (p-1)}{N+\al}(\Omega)}\le C$$
uniformly in $n$. Finally, since the embedding is compact (refer to Theorem \ref{teo appendix 1}) we have that up to subsequences 
$$\|u_n-u_\al\|_{L^\frac{2Nr}{N+\al}(\Omega)}\to 0.$$
This proves that 
$$\Big\| |u_n|^p-|u_\al|^p\Big\|_{L^\frac{2N}{N+\al}(\Omega)}\to0$$
and consequently that
$$\int\limits_\Omega(I_\al\ast |u_\al|^p)|u_\al|^p dx=1.$$
We can conclude that $u_\al\in M_\al$. Moreover, eventually passing to $|u_\al|$, we may assume this minimizer is non-negative. 
\endproof

\begin{prop}\label{esist minimo pallaC neumann}
Under the assumptions of Proposition \ref{esist minimo anello neumann}, let $\Omega=\{x\in\R^N\,|\,|x|>a>0\}$. Then, for every $p\in\left(\frac{N+\al}{N},+\infty\right)$ there exists $u_\al\in M_\al$ non-negative function such that 
$$Q(u_\al)=J_\al:=\inf\limits_{u\in M_\al} Q(u).$$
\end{prop}

\proof 
We proceed in the same way as the previous proof, but since $\Omega=\{x\in\R^N\,|\,|x|>a>0\}$, we use the following Sobolev Embedding (see Theorem \ref{teo appendix 2}): for every $q\in[2,+\infty)$ we have the continuous immersion 
$$H^1_{rad}(\Omega)\hookrightarrow L^q(\Omega)$$
which is also compact for $q\in(2,+\infty)$. Since $\{u_n\}$ is bounded in $H^1_{rad}(\Omega)$, it follows that for every $p\in\left(\frac{N+\al}{N},+\infty\right)$
$$\|u_n\|^p_ {L^\frac{2Np}{N+\al}(\Omega)}+\|u_\al\|^p_{L^\frac{2Np}{N+\al}(\Omega)}\le C$$
uniformly in $n$. Recalling estimate \eqref{formula stima diff ap-bp}, in this case we have to require that
$$\frac{2Nr'(p-1)}{N+\al}\ge2\qquad \text{and}\qquad \frac{2Nr}{N+\al}>2.$$
If $p\ge2$, taking $r=r'=2$, the two previous conditions are satisfied. If
$\frac{N+\al}{N}<p<2$, we set $r'=\frac{N+\al}{N(p-1)}$ and consequently $r=\frac{N+\al}{2N+\al-Np}$, in this way 
$$\frac{2Nr'(p-1)}{N+\al}=2\qquad\text{and}\qquad \frac{2Nr}{N+\al}=\frac{2N}{2N+\al-Np}>2$$
which concludes the proof.
\endproof

\proof[Proof of Theorem \ref{teo esist Neumann}]
Let $u_\al\in M_\al$ be a minimizer for $J_\al$ (see Proposition \ref{esist minimo anello neumann} and Proposition \ref{esist minimo pallaC neumann}). It is classical that
$u_\al$ solves
$$\begin{cases}
-\Delta u_\al+V(x)u_\al=\mu_\al(I_\al\ast|u_\al|^p)|u_\al|^{p-2}u_\al\quad \mathrm{in}\,\,\,\Omega\\
u_\al>0\hspace{6cm}\mathrm{in}\,\,\,\Omega\\
\frac{\partial u_\al}{\partial\nu}=0\hspace{5.05cm}\qquad\,\mathrm{on}\,\,\, \partial\Omega
\end{cases}$$
where $\mu_\al$ is a Lagrange multiplier, and $u_\al>0$ by the Strong Maximum Principle.  Therefore $|\mu_\al|^\frac{1}{2p-2}\,u_\al$ provides a solution to \eqref{choquad neumann}.
\endproof

\proof[Proof of Theorem \ref{teo esist Dirichlet}] When dealing with homogeneous Dirichlet boundary conditions we will investigate existence of minima of the functional $Q$ over the constrained set 
$$M_{\al,0}:=\bigg\{u\in H^1_{0,rad}(\Omega)\,\,\bigg|\,\,\int_{\Omega}\left(I_\al\ast |u|^p\right)|u|^p=1\bigg\}.$$
Exploiting Corollary \ref{cor emb H_0 anello} and Corollary \ref{cor emb H_0 fuoripalla}, one can prove that there exists $u_{\al,0}$ non-negative function which achieves the infimum 
$$J_{\al,0}:=\inf\{Q(u)\,|\, u\in M_{\al,0}\}.$$
Hence, the results of Proposition \ref{esist minimo anello neumann} and Proposition \ref{esist minimo pallaC neumann} hold also for Dirichlet boundary conditions, and we can conclude by rescaling. 
Let us assume now that $V\ge0$, $N\ge3$ and $\Omega=A_{a,b}$. If $\{u_n\}$ is a minimizing sequence for $J_{\al,0}$, exploiting estimate \eqref{stima solo grad anello} we get that
$$\int_{A_{a,b}} u_n^2\,dx\le C\int_{A_{a,b}}|\nabla u_n|^2 dx\le C\,Q(u_n)< C(J_{\al,0}+1).$$
Hence, if $N\ge3$ we obtain existence of solutions to the problem \eqref{choquad dirichlet} on $\Omega=A_{a,b}$ also in the more general case $V(|x|)\ge0$.
\endproof

%--------------------------------------------------------------------------------

\section{Nonexistence result}

First, we prove a suitable version of the celebrated Pohozaev identity (refer to the seminal papers \cite{PucciSerrin} and also \cite{GaoYang} for the Choquard case).

\begin{theorem}
Let $\Omega$ be a  smooth domain in $\R^N$. Assume that $u\in H^1_0(\Omega)\,\cap W^{2,2}(\Omega)\,\cap W^{1,\frac{2Np}{N+\al}}(\Omega)$ is a solution to 
\begin{equation}\label{Choquard p non esist}
-\Delta u+V(x)u=(I_\al\ast|u|^p)u^{p-2}u,\quad\text{in}\,\,\,\Omega
\end{equation}
such that $u^2 V\in L^1(\Omega)$ and $u^2\nabla V\cdot x\in L^1(\Omega)$. Then the following Pohozaev-type identity holds
\begin{equation}\label{pohozaev limitato}
\left(2-N+\frac{\al+N}{p}\right)\int_\Omega|\nabla u|^2dx-\left(N-\frac{\al+N}{p}\right)\int_{\Omega}Vu^2\,dx-\int_{\Omega}u^2\nabla V\cdot x\,dx=I_{\partial\Omega},
\end{equation}
where $\nu$ is the exterior unit normal at $\partial\Omega$ and $u_\nu:=\frac{\partial u}{\partial\nu}$ and $I_{\partial\Omega}=\int_{\partial\Omega}u_\nu^2\,(x\cdot\nu)\,d\s$.
\end{theorem}

\proof
Since $u$ solves \eqref{Choquard p non esist}, multiplying each term by $\nabla u\cdot x$ and integrating over $\Omega$, we get
\begin{equation}\label{14}
-\int_{\Omega}\Delta u\,\nabla u\cdot x\,dx+\int_{\Omega}Vu\,\nabla u\cdot x\,dx= \int_{\Omega}(I_\al\ast |u|^p)|u|^{p-2}u\,\nabla u\cdot x\,dx.
\end{equation}
We take into account each term of \eqref{14} separately. Integrating by parts the first term, we have
\begin{align}\notag
-\int_{\Omega}\Delta u\,\nabla u\cdot x\,dx&=\int_{\Omega}\nabla u\cdot\nabla(\nabla u\cdot x)dx-\int_{\partial\Omega}(\nabla u\cdot \nu)(\nabla u\cdot x)\,d\s\\ \label{primo pezzo poho}
&=\int_{\Omega}|\nabla u|^2+\int_{\Omega}\sum \limits_{i,j}u_{x_i}u_{x_i\,x_j}x_j-\int_ {\partial\Omega}(\nabla u\cdot \nu)(\nabla u\cdot x)\, d\s.
\end{align}
Notice that integrating by parts 
\begin{equation}\label{identita dim poho}
\int_{\Omega}\sum\limits_{i,j}(u_{x_i}x_j)u_{x_ix_j}=\int_{\partial\Omega}|\nabla u|^2\,x\cdot\nu\,d\s-\int_{\Omega}\nabla u\cdot\nabla(\nabla u\cdot x)+(1-N)\int_ {\Omega}|\nabla u|^2.
\end{equation}
Substituting \eqref{identita dim poho} into \eqref{primo pezzo poho}, and integrating by parts again, we find
\begin{align*}
-\int_{\Omega}\Delta u\,\nabla u\cdot x\,dx=(2-N)&\int_{\Omega}|\nabla u|^2\,dx +\int_{\Omega}\Delta u\,\nabla u\cdot x\,dx+\\
&+\int_{\partial\Omega}|\nabla u|^2\,x\cdot\nu\,d\s-2\int_ {\partial\Omega}(\nabla u\cdot \nu)(\nabla u\cdot x)\, d\s.
\end{align*}
Since $u=0$ on $\partial\Omega$, one has that $\nabla u(x)=u_\nu\nu$ where $u_\nu=\frac{\partial u}{\partial\nu}$, hence we get
\begin{equation}\label{1 poho2}
-\int_{\Omega}\Delta u\,\nabla u\cdot x\,dx=\frac{2-N}{2}\int_{\Omega}|\nabla u|^2\,dx -\frac{1}{2}\int_{\partial\Omega}u_\nu^2\,x\cdot\nu\,d\s.
\end{equation}
Concerning the second term, integrating by parts we have
$$\int_{\Omega}Vu\,\nabla u\cdot x\,dx=\int_{\partial\Omega}Vu^2\,x \cdot\nu\,d\s-\int_{\Omega}Vu\,\nabla u\cdot x\,dx-N\int_{\Omega}Vu^2\,dx-\int_{\Omega}u^2\nabla V\cdot x\,dx$$
and hence
\begin{equation}\label{2 poho2}
\int_{\Omega}Vu\,\nabla u\cdot x\,dx=\frac{1}{2}\int_{\partial\Omega}Vu^2\,x\cdot\nu\, d\s-\frac{N}{2}\int_{\Omega}Vu^2\,dx-\frac{1}{2}\int_{\Omega}u^2\nabla V\cdot x\,dx.
\end{equation}
As for the Riesz term we get
\begin{align}\notag
\int\limits_{\Omega}(I_\al&\ast|u|^p)|u|^{p-2}u\,\nabla u\cdot x\,dx= c\int\limits_{\Omega}\int \limits_{\Omega}\frac{|u(y)|^p\, |u(x)|^{p-2}u(x)\nabla u\cdot x}{|x-y|^{N-\al}}dy\,dx=\\\label{cinque 2}
&=c\int\limits_{\Omega}\int\limits_{\partial\Omega}\frac{|u(y)|^p\,|u(x)|^p}{|x-y|^{N-\al}}(x\cdot\nu)\,d\s(x)\,dy-c\int\limits_{\Omega}\int\limits_{\Omega} |u(y)|^p\,u(x)\,\mathrm{div}_x \left(\frac{u(x)|u(x)|^{p-2}x}{|x-y|^{N-\al}}\right)dx\,dy,
\end{align}
where $c=c(N,\al)$ is the constant in the definition of the Riesz potential. We have that
\begin{align}\notag
\int\limits_{\Omega}\int\limits_{\Omega}&|u(y)|^p u(x)\,\mathrm{div}_x \left(\frac{u(x)|u(x)|^{p-2}}{|x-y|^{N-\al}}x\right)dx\,dy\\ \notag
&=\int\limits_{\Omega} \int\limits_\Omega\frac{|u(y)|^p\,u(x)|u(x)|^{p-2}\nabla u\cdot x}{|x-y|^{N-\al}}dx\,dy+(p-2)\int\limits_{\Omega}\int\limits_ \Omega\frac{|u(y)|^p\,u(x)|u(x)|^{p-2}\nabla u\cdot x}{|x-y|^{N-\al}}dx\,dy\\ \notag
&\quad+(\al-N)\int\limits_{\Omega}\int\limits_{\Omega}\frac{|u(x)|^p\,|u(y)|^p}{|x-y|^{N-\al}}\frac{(x-y)\cdot x}{|x-y|^2}dx\,dy+N\int\limits_{\Omega}\int\limits_ {\Omega}\frac{|u(x)|^p\,|u(y)|^p}{|x-y|^{N-\al}}dx\,dy\\ \notag
&=(p-1)\int\limits_{\Omega}\int\limits_\Omega\frac{|u(y)|^p\, u(x)|u(x)|^{p-2}\nabla u\cdot x}{|x-y|^{N-\al}}dx\,dy+\frac{\al+N}{2}\int\limits_{\Omega}\int\limits_ {\Omega}\frac{|u(x)|^p\,|u(y)|^p}{|x-y|^{N-\al}}dx\,dy\\\label{divergenza riesz}
&\quad+\frac{\al-N}{2}\int\limits_{\Omega}\int\limits_{\Omega}\frac{|u(x)|^p\, |u(y)|^p}{|x-y|^{N-\al}}\frac{(x+y)\cdot (x-y)}{|x-y|^2}dx\,dy
\end{align}
where we used that $\frac{x\cdot(x-y)}{|x-y|^2}=\frac{1}{2}+\frac{(x+y) \cdot(x-y)}{2|x-y|^2}$ and moreover we observe that by symmetry 
$$\int\limits_{\Omega}\int\limits_{\Omega}\frac{|u(x)|^p\,|u(y)|^p}{|x-y|^{N-\al}}\frac{(x+y)\cdot (x-y)}{|x-y|^2}dx\,dy=0.$$
Using \eqref{divergenza riesz} in \eqref{cinque 2} we finally get
\begin{align}\notag
\int\limits_{\Omega}(I_\al\ast|u|^p)&|u|^{p-2}u\,\nabla u\cdot x\,dx\\ \label{3 poho2}
&=\frac{c}{p}\int\limits_{\Omega}\int\limits_{\partial\Omega}\frac{|u(y)|^p\,|u(x)|^p}{|x-y|^{N-\al}}(x\cdot\nu)\,d\s(x)\,dy-\frac{\al+N}{2p}\int\limits_{\Omega}(I_\al\ast|u|^p)|u|^p dx.
\end{align}
Summing up \eqref{1 poho2}, \eqref{2 poho2} and \eqref{3 poho2}, and using the fact that $u=0$ on $\partial\Omega$ we obtain the following identity
\begin{equation}\label{poho2 lunga}
\begin{split}
 \frac{2-N}{2}\int_{\Omega}|\nabla u|^2\,dx-\frac{N}{2}\int_{\Omega}&Vu^2\,dx-\frac{1}{2}\int_{\Omega}u^2\nabla V\cdot x\,dx\\
 &+\frac{\al+N}{2p}\int\limits_{\Omega}(I_\al\ast|u|^p)|u|^p\,dx=\frac{1}{2}\int_{\partial\Omega}u_\nu^2\,x\cdot\nu\,d\s.   
\end{split}
\end{equation}
Finally, testing equation \eqref{Choquard p non esist} with $u$, we infer that
$$\int_\Omega(I_\al\ast |u|^p)|u|^p\,dx=\int_\Omega|\nabla u|^2dx+\int_\Omega Vu^2\,dx$$
using this relation in \eqref{poho2 lunga} we conclude the proof of the Pohozaev-type identity \eqref{pohozaev limitato}.
\endproof

\begin{obs}
Notice that the previous integrability assumptions on $u$ can be weakened for suitable values of the nonlinearity $p$, indeed one can adapt the proof of regularity for solutions to the Choquard equation in the whole space $\R^N$ to the case of  domains (see \cite[Proposition 4.1]{MVS13} and also \cite[Theorem 2]{MVS15}).
\end{obs}

As a consequence, using the Pohozaev identity we can prove that if the domain $\Omega$ is strictly star-shaped, the value $p=\frac{N+\al}{N-2}$ is critical from the point of view of existence of non-trivial solution for the Dirichlet problem \eqref{pb non exist}. Notice that in the literature also $p=\frac{\al+N}{N}$ is a critical value for the Choquard equation in the whole space $\R^N$ and for Choquard boundary value problems defined in exterior domains of the form $\R^N\setminus \Bar{B_a}(0)$ (refer to Theorem \ref{teo esist Neumann} and Theorem \ref{teo esist Dirichlet}). 

\proof[Proof of Theorem \ref{teo non esist}]
Since $u\in H^1_0(\Omega)\,\cap W^{2,2}(\Omega)\,\cap W^{1,\frac{2Np}{N+\al}}(\Omega)$ 
solves \eqref{pb non exist} and we assumed that $u^2 V, \,\,u^2\nabla V\cdot x\in L^1(\Omega)$, we have the following Pohozaev identity
$$\left(N-2-\frac{\al+N}{p}\right)\int_\Omega|\nabla u|^2dx+\left(N-\frac{\al+N}{p}\right)\int_{\Omega}Vu^2\,dx+\int_{\Omega}u^2\nabla V\cdot x\,dx+\int_{\partial \Omega}u_\nu^2\,(x\cdot\nu)\,d\s=0$$
Since $\Omega$ is strictly star-shaped with respect to $0\in\R^N$, we have that $x\cdot\nu>0$ on $\partial\Omega$, moreover by assumptions $\int_{\Omega}u^2\nabla V\cdot x\,dx\ge0$. So if $V\geq 0$ and  $p>\frac{N+\al}{N-2}$, $u$ must be identically equal to $0$, and in the case $p=\frac{N+\al}{N-2}$ the same holds if $V>0$. 
\endproof

%------------------------------------------------------------------------------------

\section{Limiting problem}

In this section we will always assume that $\Omega=A_{a,b}$, $N\ge2$, $p$ is a fixed real value in $[1,+\infty)$ and $V(|x|)$ is a continuous radial function on $\Omega$ such that $\inf_{x\in\Omega}V(x)>0$. Notice that if $0\le\al_2<\al_1<N$ for every couple of points $x,y\in B_b(0)$ it holds
$$\frac{1}{|x-y|^{N-\al_1}}\le(\max\{1,2b\})^{\al_1-\al_2}\frac{1}{|x-y|^{N-\al_2}},$$
from which using that $C_{N,\al_1}<C_{N,\al_2}$ (as defined in \eqref{def riesz potential}) it follows that
\begin{equation}\label{stima alfa1>alfa2}
\int_{\Omega}(I_{\al_1}\ast |f|)|f|\,dx\le C\int_{\Omega}(I_{\al_2}\ast |f|)|f|\,dx
\end{equation}
where $C=(\max\{1,2b\})^{\al_1}$. From Proposition \ref{esist minimo anello neumann} we have that for every fixed $\al\in(0,N)$, there exists a non-negative function $u_\al\in M_\al$ which minimizes $Q(u)$, that is $J_\al=Q(u_\al)$. 

We recall the following classical result (see \cite[Theorem D]{kuro}). 

\begin{theorem}[Riesz kernel as an approximation of the identity]\label{riesz alpha 0 identità}]   Let $f\in L^p(\R^N)$ for $p\in[1,+\infty)$. If the Riesz potential $I_\al\ast f$ is well-defined, we have that 
$$I_\al\ast f(x)\to f(x), \quad\text{as}\,\,\,\al\to0^+$$
at each Lebesgue point of $f$.    
\end{theorem}

Following some arguments of \cite{GrossiNoris}, we prove weak convergence of minimizers as $\al\to0^+$. 

\begin{prop}\label{prop 4.1}
Let $\{\al_n\}$ be a sequence of real values in $(0,N)$ such that $\al_n\to0^+$ as $n\to+\infty$ and $u_{\al_n}\in M_{\al_n}$ be a sequence of minimizers to $J_{\al_n}$ respectively. Then, there exists $u_0\in H_{rad}^1(\Omega)$ such that, up to subsequences 
$$u_{\al_n}\rightharpoonup u_0\quad\text{in}\,\,\,H^1_{rad}(\Omega)\qquad u_{\al_n}\to u_0\quad\text{in}\,\,\,L^q(\Omega),\quad\forall q\in[1,+\infty)$$
 and \[ \int\limits_{\Omega}u_0^{2p}(x)\,dx=1.\] 
 \end{prop}

\proof
First we prove that the sequence $\{u_{\al_n}\,|\, n\in\mathbb{N}\,\}$ is bounded in $H^1_{rad}(\Omega)$. Let us consider a non-negative test function $\eta\in M_{\al_1}$ and for every $\al\in(0,\al_1)$ let us define
$$\eta_\al:=\frac{\eta}{\left(\int\limits_\Omega(I_\al\ast \eta^p)\eta^p\,dx\right)^\frac{1}{2p}}.$$
Denoting by $a=\Big(\int\limits_\Omega(I_\al\ast\eta^p)\eta^p\,dx\Big)^\frac{1}{2p}$ we observe that
$$\int\limits_\Omega (I_\al\ast\eta_\al^p)\eta_\al^p\,dx=\frac{1}{a^{2p}}\int\limits _\Omega(I_\al\ast\eta^p)\eta^p\,dx=1$$
hence $\eta_\al\in M_\al$. For every $\al\in(0,\al_1)$, using \eqref{stima alfa1>alfa2}, we get
$$Q(\eta_\al)=\frac{Q(\eta)}{a^2}=\frac{Q(\eta)}{\Big(\int\limits_\Omega(I_\al\ast\eta ^p)\eta^p\,dx\Big)^\frac{1}{p}}\le\frac{Q(\eta)}{\Big(\frac{1}{C}\int\limits _\Omega(I_{\al_1}\ast\eta^p)\eta^p\,dx\Big)^\frac{1}{p}}=C_1\,Q(\eta)$$
where $C_1=C_1(b,\al_1,p)$ and in the last equality we exploited the fact that $\eta\in M_{\al_1}$. Since $\eta_\al\in M_\al$, it follows that $J_\al\le C_1\,Q(\eta)$, therefore $J_\al\le C_1$ for every $\al\in(0,\al_1)$. Hence if $u_{\al_n}$ is a minimizer associated to $\al_n$, we have that
$$\int_{\Omega}|\nabla u_{\al_n}|^2\le C_1\qquad\text{and}\qquad\int_{\Omega}|u_ {\al_n}|^2\le C_1$$
uniformly in $n$. This proves that the sequence $\{u_{\al_n}\,|\,n\in\mathbb{N}\}$ is bounded in $H^1_{rad}(\Omega)$ as $\al_n\to0^+$, so there exists $u_0\in H_{rad}^1(\Omega)$ such that, up to subsequences $u_{\al_n}\rightharpoonup u_0$ in $H^1_{rad}(\Omega)$. We can conclude using the compact Sobolev Embedding for radial functions (refer to Theorem \ref{teo appendix 1}).

Now  we prove that $\int\limits_\Omega u_0^{2p}(x)\,dx\ge1$. We observe that if $u_{\al_1}\in M_{\al_1}$, using \eqref{stima alfa1>alfa2} it holds
\begin{equation}\label{stima 1<=al2 to 0}
1=\int\limits_{\Omega}(I_{\al_1}\ast|u_{\al_1}|^p)|u_{\al_1}|^p\,dx\le\lim\limits_ {\al_2 \to0^+}C\int\limits_\Omega(I_{\al_2}\ast|u_{\al_1}|^p)|u_{\al_1}|^p\,dx.
\end{equation}
Passing to the limit as $\al_1\to0^+$ in \eqref{stima 1<=al2 to 0}, we get that
$$1=\lim\limits_{\al_1\to 0^+}\int_{\Omega}(I_{\al_1}\ast|u_{\al_1}|^p)|u_{\al_1}| ^p\,dx\le\lim\limits_{\al_1\to0^+}\lim\limits_{\al_2\to0^+}C\int\limits_\Omega(I_{\al_2} \ast|u_{\al_1}|^p)|u_{\al_1}|^p\,dx,$$
and using monotone convergence and Theorem \ref{riesz alpha 0 identità}, we obtain that
$$1\le\lim\limits_{\al_1\to0^+}C\int\limits_\Omega\lim\limits_{\al_2\to0^+}(I_{\al_2} \ast|u_{\al_1}|^p)|u_{\al_1}|^p\,dx=\lim\limits_{\al_1\to0^+}C\int\limits_\Omega |u_{\al_1}(x)|^{2p}\,dx.$$
Since $u_{\al_n}\to u_0$ in $L^q(\Omega)$ for every $q\in[1,+\infty)$ and $C=(\max\{1,2b\})^{\al_1}\to1$ as $\al_1\to0$, we can conclude that 
$$1\le\int\limits_\Omega u_0(x)^{2p}\,dx.$$
As for the inverse inequality, we observe that again by \eqref{stima alfa1>alfa2}, for $\al_2>\al_1$
$$\lim\limits_{\al_1\to 0^+}\int\limits_\Omega(I_{\al_2}\ast|u_{\al_1}|^p)|u_{\al_1}| ^p dx\le\lim\limits_{\al_1\to0^+}C\int\limits_\Omega(I_{\al_1}\ast|u_{\al_1}|^p) |u_{\al_1}| ^p dx=C$$
where $C=(\max\{1,2b\})^{\al_2}$. Passing to the limit as $\al_2\to0^+$ we get
$$1\ge\lim\limits_{\al_2\to 0^+}\lim\limits_{\al_1\to0^+}\int\limits_\Omega(I_{\al_2} \ast|u_{\al_1}|^p)|u_{\al_1}|^p dx=\lim\limits_{\al_2\to0^+}\int\limits_ \Omega(I_{\al_2}\ast|u_0|^p)|u_0|^p dx$$
where to get the last equality, we proceed as in the proof of Proposition \ref{esist minimo anello neumann}. Finally, using again the Monotone Convergence Theorem and Theorem \ref{riesz alpha 0 identità}, we obtain that
$$1\ge\int\limits_\Omega u_0(x)^{2p}\,dx$$
which concludes the proof.
\endproof

Now we prove that every non-negative function belonging to $M_0$ can be seen as the limit of a suitable approximating sequence.

\begin{prop}\label{prop 4.3}
Let $u\in M_0$ be a non-negative function. Then for every $\al\in(0,N)$ there exists $w_\al\in M_\al$ such that $w_\al\to u$ in $H^1_{rad}(\Omega)$ as $\al\to 0^+$.
\end{prop}

\proof
For every $\al\in(0,N)$, let us consider the function
$$w_{\s_\al}:=\s_\al u,\qquad\text{where}\quad\s_\al=\left(\int_\Omega(I_\al\ast|u|^p)|u|^p \right)^{-\frac{1}{2p}}.$$
It follows immediately that
$$\int\limits_\Omega(I_\al\ast |w_{\s_\al}|^p)|w_{\s_\al}|^p dx=\s_\al^{2p}\int \limits_\Omega (I_\al\ast |u|^p)|u|^p dx=1,$$
hence $w_{\sigma_\al}$ belongs to $M_\al$.  In order to conclude we observe that
$$\lim\limits_{\al\to0^+}\s_\al=\lim\limits_{\al\to0^+}\left(\int_\Omega(I_\al \ast |u|^p)|u|^p dx\right)^{-\frac{1}{2p}}=\left(\int_\Omega |u(x)|^{2p}\,dx\right)^{-\frac{1}{2p}}=1$$
using in the last equality that $u\in M_0$.
\endproof

In order to conclude we have to show that the constrained variational problems 
$J_\al$ converges to the limit problem $J_0$ as $\al\to0^+$.

\begin{prop}
Let us consider a sequence of real values $\al\in(0,N)$ such that $\al\to0^+$ and let $\{u_{\al}\}$ be a sequence of minimizers of $Q(u)$ in $M_{\al}$ respectively. Then, up to subsequences, we have that 
$$u_{\al}\to u_0\quad\,\text{in} \,\,\,H^1_{rad}(\Omega)$$
and 
$$J_{\al}=\inf\limits_{u\in M_{\al}}Q(u)=Q(u_{\al})\xrightarrow[\al\to0^+]{}J_0=\inf \limits_{u\in M_0}Q(u)=Q(u_0).$$
\end{prop}

\proof 
From Proposition \ref{prop 4.1}  there exists $u_0\in M_0$ such that $u_{\al}\rightharpoonup u_0$ in $H^1_{rad}(\Omega)$. So, by lower semicontinuity with respect to weak convergence, we get that
\begin{equation}\label{J0<liminf}
J_0\le Q(u_0)\le\liminf\limits_{\al\to0^+}Q(u_\al)=\liminf\limits_{\al\to0^+}J_\al \end{equation}
where we used also that $J_0:=\inf\limits_{u\in M_0} Q(u)$ and that every $u_\al$ is a minimizer for $J_\al$. 
Conversely, let $u\in M_0$ be a non-negative function, by Proposition \ref{prop 4.3} there exists a sequence of function $w_\al\in M_{\al}$ such that $w_\al\to u$ in $H^1_{rad}(\Omega)$ as $\al\to 0^+$, so it holds
\begin{equation}\label{Q>limsup}
Q(u)=\lim\limits_{\al\to0^+} Q(w_\al)\ge\limsup\limits_{\al\to0^+}J_\al.    
\end{equation}
Notice that in the minimization problem $J_0=\inf\limits_{u\in M_0}Q(u)$, we can equivalently minimize over the constrained set $\{u \in M_0\,\,|\,\, u\ge0\}$ (possibly taking $|u|$ instead of $u$). Passing to the infimum in \eqref{Q>limsup}, we get that
$$J_0=\inf\limits_{\{u \in M_0|u\ge0\}} Q(u)\ge\limsup\limits_ {\al\to0^+}J_\al$$
and hence from \eqref{J0<liminf} we finally obtain
$$J_0\ge\limsup\limits_ {\al\to0^+}J_\al\ge\liminf\limits_{\al\to0^+}J_\al\ge J_0,$$
which concludes the proof.
\endproof

\proof[Proof of Theorem \ref{teo problema limite}]
Since $u_0$ is a minimizer for $J_0$, we have that
$$\begin{cases}
-\Delta u_0+V(x)u_0=\mu_0 u_0^{2p-1}\,\qquad\quad \mathrm{in}\,\,\,\Omega\\
u_0>0\hspace{4.28cm}\mathrm{in}\,\,\,\Omega\\
\frac{\partial u_0}{\partial\nu}=0\hspace{3.3cm}\qquad\,\mathrm{on}\,\,\, \partial\Omega
\end{cases}$$
where $\mu_0$ is a Lagrange multiplier. A suitable multiple of $u_0$ (namely $|\mu_0|^\frac{1}{2p-2}\,u_0$) provides a solution to \eqref{pb u2p}.
\endproof

%------------------------------------------------------------------------------------


\begin{thebibliography}{90}

\bibitem{Adams} Adams R.A.; \textit{Sobolev spaces.} Pure and Applied Mathematics, Vol. 65. Academic Press [Harcourt Brace Jovanovich, Publishers], New York-London, 1975. xviii+268 pp.

\bibitem{AmbrosettiMalchiodi}  Ambrosetti A., Malchiodi A.; \textit{Nonlinear analysis and semilinear elliptic problems}. Cambridge Studies in Advanced Mathematics, 104. Cambridge University Press, Cambridge, 2007. xii+316 pp. ISBN: 978-0-521-86320-9; 0-521-86320-1.

\bibitem{Bernardini4} Bernardini C.: \textit{Mass concentration for Ergodic Choquard Mean-Field Games.} (2022) submitted (preprint ArXiv:\,\href{https://arxiv.org/abs/2212.00132}{\texttt{2212.00132}})

\bibitem{BerCes} Bernardini C., Cesaroni A.; \textit{Ergodic Mean-Field Games with aggregation of Choquard-type,} J. Differential Equations \textbf{364}, 296-335 (2023) \doi{10.1016/j.jde.2023.03.045}.

\bibitem{Brezis}  Brezis H.; \textit{Functional analysis, Sobolev spaces and partial differential equations}. Universitext. Springer, New York, 2011. xiv+599 pp. ISBN: 978-0-387-70913-0

\bibitem{ChenLiu}  Chen P., Liu X.; \textit{Positive solutions for Choquard equation in exterior domains}. Commun. Pure Appl. Anal. 20 (2021), no. 6, 2237–2256. 

\bibitem{ClappSal} Clapp M., Salazar D.; \textit{Positive and sign changing solutions to a nonlinear Choquard equation.} J. Math. Anal. Appl. 407 (2013), no. 1, 1–15. 

\bibitem{GaoYang} Gao F., Yang M.; \textit{The Brezis-Nirenberg type critical problem for the nonlinear Choquard equation}, Sci China Math, 2018, 61: 1219–1242, \doi{10.1007/s11425-016-9067-5}.

\bibitem{GhiPag} Ghimenti M., Pagliardini D.; \textit{Multiple positive solutions for a slightly subcritical Choquard problem on bounded domains}, Calc. Var. Partial Differential Equations 58 (2019), no. 5, Paper No. 167, 21 pp.

\bibitem{Goel1} Goel D., Rădulescu V.D., Sreenadh, K.; \textit{Coron problem for nonlocal equations involving Choquard nonlinearity.} Adv. Nonlinear Stud. 20 (2020), no. 1, 141–161.

\bibitem{Goel2} Goel D., Sreenadh K.; \textit{Critical growth elliptic problems involving Hardy-Littlewood-Sobolev critical exponent in non-contractible domains}, Adv. Nonlinear Anal. 9 (2020), no. 1, 803–835.

\bibitem{Grossi} Grossi M.; \textit{Asymptotic behaviour of the Kazdan-Warner solution in the annulus.} J. Differential Equations 223 (2006), no. 1, 96–111. 

\bibitem{GrossiNoris} Grossi M., Noris B.; \textit{Positive constrained minimizers for supercritical problems in the ball}, Proc. Amer. Math. Soc., 140:2141–2154, (2012).

\bibitem{KadzanWarner} Kazdan J.L., Warner F.W.; \textit{Remarks on some quasilinear elliptic equations}, Comm. Pure Appl.Math. 28 (1975) 567–597.

\bibitem{kuro} Kurokawa T.; \textit{On the riesz and bessel kernels as approximations of the identity}. Sci. Rep. Kagoshima Univ, 30:31–45, 1981.

%\bibitem{Leoni}  Leoni G.; \textit{A first course in Sobolev spaces.} Graduate Studies in Mathematics, 105. American Mathematical Society, Providence, RI, 2009. xvi+607 pp. ISBN: 978-0-8218-4768-8.

\bibitem{Lieb} Lieb E.H.; \textit{Existence and uniqueness of the minimizing solution of Choquard's nonlinear equation}, Studies in Appl. Math. 57 (1976/77), no. 2, 93–105.

%\bibitem{LL} Lieb E.H., Loss M.; \textit{Analysis.} Second edition. Graduate Studies in Mathematics, 14. American Mathematical Society, Providence, RI, 2001. xxii+346 pp. ISBN: 0-8218-2783-9 00A05 (26-01 28-01 31-01 35J10 42-01)

\bibitem{Lions} Lions P.L.; \textit{The Choquard equation and related questions}, Nonlinear Anal. 4 (1980), no. 6, 1063–1072.

%\bibitem{Lions82} Lions P.L.; \textit{Compactness and topological methods for some nonlinear variational problems of mathematical physics}, Nonlinear problems: present and future (Los Alamos, N.M., 1981), North-Holland Math. Stud., vol. 61, North-Holland, Amsterdam–New York, 1982, pp. 17–34. 

%\bibitem{Lions84} Lions P.L.; \textit{The concentration-compactness principle in the calculus of variations. The locally compact case. I.} Ann. Inst. H. Poincaré Anal. Non Linéaire, 1(2):109–145, 1984.

\bibitem{MaZhao} Ma L., Zhao L.; \textit{Classification of positive solitary solutions of the nonlinear Choquard equation}, Arch. Ration. Mech. Anal. 195 (2010), no. 2, 455–467, \doi{10.1007/s00205-008-0208-3}.

%\bibitem{Menzala} Menzala G.P.; \textit{On regular solutions of a nonlinear equation of Choquard's type}, Proc. Roy. Soc. Edinburgh Sect. A 86 (1980), no. 3–4, 291–301, \doi{ 10.1017/S0308210500012191}. 

\bibitem{MVS13} Moroz V., Van Schaftingen J.; \textit{Groundstates of nonlinear Choquard equations: existence, qualitative properties and decay asymptotics}, J. Funct. Anal. 265 (2013), no. 2, 153–184, \doi{10.1016/j.jfa.2013.04.007}. 

\bibitem{MVS13bis} Moroz V., Van Schaftingen J.; \textit{Nonexistence and optimal decay of supersolutions to Choquard equations in exterior domains}, J. Diff. Eq. 254, Vol.8,
2013. \doi{10.1016/j.jde.2012.12.019}.

\bibitem{MVS15} Moroz V., Van Schaftingen J.; \textit{Existence of groundstates for a class of nonlinear Choquard equations}, Trans. Am. Math. Soc. 367(9), 6557–6579 (2015)

\bibitem{MVS15 2} Moroz V.; Van Schaftingen J.; \textit{Groundstates of nonlinear Choquard equations: Hardy-Littlewood-Sobolev critical exponent}. Commun. Contemp. Math. 17 (2015), no. 5, 1550005, 12 pp.

\bibitem{MS} Moroz V., Van Schaftingen J.; \textit{A guide to the Choquard equation}, J. Fixed Point Theory Appl. \textbf{19} (2019), 773–813. \doi{10.1007/s11784-016-0373-1}

\bibitem{NiNuss} Ni W.M., Nussbaum R.; \textit{Uniqueness and nonuniqueness for positive radial solutions of $\Delta u+f(u, r)=0$}, Comm. Pure Appl. Math. 38 (1985) 67–108.

\bibitem{Pekar} Pekar S.; \textit{Untersuchung über die Elektronentheorie der Kristalle}, Akademie Verlag, Berlin, 1954.

%\bibitem{Pohozaev} Pohožaev S.I.; \textit{On the eigenfunctions of the equation $\Delta u+\la f(u)=0$}, Dokl. Akad. Nauk SSSR 165, 36–39 (1965) (Russian); English transl., Soviet Math. Dokl. 6, 1408–1411 (1965).

\bibitem{PucciSerrin} Pucci P., Serrin J.; \textit{A general variational identity}, Indiana Univ. Math. J. 35(3), 681–703 (1986).

\bibitem{Strauss} Strauss W.A.; \textit{Existence of solitary waves in higher dimensions.} Comm. Math. Phys. 55 (1977), no. 2, 149–162. 

\bibitem{Struwe} Struwe M.; \textit{Variational methods. Applications to nonlinear partial differential equations and Hamiltonian systems}. Second edition. Ergebnisse der Mathematik und ihrer Grenzgebiete (3), 34. Springer-Verlag, Berlin, 1996. xvi+272 pp. 
ISBN: 3-540-58859-0

\bibitem{VanSchXia} Van Schaftingen J., Xia J.; \textit{Choquard equations under confining external potentials}, Nonlinear Differential Equations Appl. 24 (2017), no. 1, Paper No. 1, 24 pp.

%\bibitem{Yang} Yang X.; \textit{Existence of positive solution for the Choquard equation in exterior domain}, Complex Variables and Elliptic Equations, 67:8, 2043-2059, \doi{10.1080/17476933.2021.1913133}

%\bibitem{WeiWinter} Wei J., Winter M.; \textit{Strongly interacting bumps for the Schrödinger–Newton equations}, J. Math. Phys. 50, 012905 (2009); \doi{10.1063/1.3060169}

\end{thebibliography}
\end{document}